\documentclass[12pt]{article}
\usepackage{amsmath, amsthm, amssymb,latexsym,enumerate}
\usepackage{graphicx}
\usepackage{hyperref}
\usepackage{xcolor}
\usepackage{todonotes}
\usepackage{cite}
\usepackage{dsfont}
\usepackage{graphicx, caption, subcaption}

\newtheorem{theorem}{Theorem}[section]
\newtheorem{lemma}[theorem]{Lemma}

\newtheorem{claim}[theorem]{Claim}

\numberwithin{equation}{section}

\newcommand{\vp}{\varphi}
\newcommand{\uc}{\vp_*}
\newcommand{\pcf}{\chi_{\mathrm{pcf}}}

\DeclareMathOperator{\mad}{mad}

\title{Proper conflict-free coloring of sparse graphs}

\author{Eun-Kyung Cho\thanks{
Department of Mathematics, Hankuk University of Foreign Studies, Yongin-si, Gyeonggi-do, Republic of Korea.
 \texttt{ekcho2020@gmail.com}
}
\and Ilkyoo Choi\thanks{
Department of Mathematics, Hankuk University of Foreign Studies, Yongin-si, Gyeonggi-do, Republic of Korea.
\texttt{ilkyoo@hufs.ac.kr}
}
\and Hyemin Kwon\thanks{
Department of Mathematics, Ajou University, Suwon-si, Gyeonggi-do, Republic of Korea.
\texttt{khmin1121@ajou.ac.kr}
}
\and Boram Park\thanks{
Department of Mathematics, Ajou University, Suwon-si, Gyeonggi-do, Republic of Korea.
\texttt{borampark@ajou.ac.kr}
}}
\date\today

\begin{document}
 
\maketitle
\begin{abstract}
     
 A {\it proper conflict-free $c$-coloring} of a graph is a proper $c$-coloring such that each non-isolated vertex has a color appearing exactly once on its neighborhood.  
 This notion was formally introduced by Fabrici et al.,
 who proved that planar graphs have a proper conflict-free 8-coloring and constructed a planar graph with no proper conflict-free 5-coloring. 
 Caro, Petru\v{s}evski, and \v{S}krekovski investigated this coloring concept further, and in particular studied upper bounds on the maximum average degree that guarantees a proper conflict-free $c$-coloring for $c\in\{4,5,6\}$.

 Along these lines, we completely determine the threshold on the maximum average degree of a graph $G$, denoted $\mad(G)$, that guarantees a proper conflict-free $c$-coloring for all $c$ and also provide tightness examples. Namely, for $c\geq 5$ we prove that a graph $G$ with $\mad(G)\leq \frac{4c}{c+2}$ has a proper conflict-free $c$-coloring, unless $G$ contains a $1$-subdivision of the complete graph on $c+1$ vertices. 
 When $c=4$, we show that a graph $G$ with $\mad(G)<\frac{12}{5}$ has a proper conflict-free $4$-coloring, unless $G$ contains an induced $5$-cycle. 
 In addition, we show that a planar graph with girth at least 5  has a proper conflict-free $7$-coloring.
\end{abstract}

\section{Introduction and notation} 
All graphs in this paper are finite and simple, which means no loops and no parallel edges. 
For a graph $G$, let $V(G)$ and $E(G)$ denote its vertex set and its edge set, respectively.
Given a graph $G$, the {\it square} of $G$ is obtained by adding all edges between vertices of distance at most 2 in $G$. 

Coloring the square of a graph is a frequently considered problem in graph coloring, see~\cite{2008KrKr} for a survey on graph coloring with distance constraints and~\cite{1995JeTo} for a prominent book on graph coloring. 
A lower bound on the chromatic number of the square of a graph $G$ is the maximum degree of $G$ plus one, since a vertex and all of its neighbors must receive distinct colors. 
As a relaxation, one could require only one neighbor to have a distinct color among a vertex and its neighbors. 
This relaxed coloring is named proper conflict-free coloring, which we formally define as follows. 
A {\it proper conflict-free $c$-coloring} (or a PCF $c$-coloring for short) of a graph is a proper $c$-coloring such that each non-isolated vertex has a color appearing exactly once on its neighborhood.  
The {\it proper conflict-free chromatic number} of a graph $G$, denoted $\pcf(G)$, is the minimum $c$ such that $G$ has a PCF $c$-coloring.  

The proper conflict-free chromatic number was formally introduced by Fabrici et al.~\cite{arXiv_FaLuRiSo}, where they investigated the proper conflict-free chromatic number of planar graphs and outerplanar graphs, among many other related variants of a proper conflict-free coloring.
In particular, they proved that a planar graph has a proper conflict-free 8-coloring, and constructed a planar graph with no proper conflict-free 5-coloring. 
Other graph classes such as minor-closed families and bounded layered treewidth were investigated by Liu~\cite{arXiv_Liu}, and graphs with bounded expansion were studied by Hickingbotham~\cite{arXiv_Hickingbotham} with corollaries regarding $k$-planar graphs and graphs with crossings on other surfaces. 

Caro, Petru\v sevski, and \v Skrekovski~\cite{arXiv_CaPeSk} investigated the proper conflict-free chromatic number in a variety of directions, including a Brooks-type result and results regarding sparse graphs, where spareness is measured in terms of the {\it maximum average degree}. 
(Given a graph $G$, the {\it maximum average degree} of $G$, denote $\mad(G)$, is the maximum of $\frac{2|E(H)|}{|V(H)|}$ over all non-empty subgraphs $H$ of $G$.)
In particular, they proved the following theorem regarding proper conflict-free colorings of graphs with bounded maximum average degree: 

\begin{theorem}[\cite{arXiv_CaPeSk}]\label{thm:CaPeSk}
Let $G$ be a graph.
\begin{enumerate}[(i)]
    \item If $\mad(G)<\frac{8}{3}$, then $\pcf(G)\leq 6$.
    \item If $\mad(G)<\frac{5}{2}$, then $\pcf(G)\leq 5$.
    \item If $\mad(G)<\frac{24}{11}$, then $\pcf(G)\leq 4$, unless every maximal $2$-connected subgraph of $G$ is a $5$-cycle. 
\end{enumerate}
\end{theorem}

Since every planar graph $G$ with girth at least $g$ satisfies $\mad(G)<\frac{2g}{g-2}$, Theorem~\ref{thm:CaPeSk} has the following corollaries: a planar graph with girth at least $8, 10, 24$ has a proper conflict-free $6$-, $5$-, $4$-coloring, respectively. 
They actually improve the girth constraint implied by $(i)$ and obtain that a planar graph with girth at least $7$ has a proper conflict-free $6$-coloring. 

In this paper, we not only improve upon the maximum average degree bound in Theorem~\ref{thm:CaPeSk}, but also completely determine the threshold on the maximum average degree that ensures a graph to have a proper conflict $c$-coloring for all $c$. 
It is not hard to see that a graph $G$ with $\mad(G)\leq c-1$ has a proper conflict-free $c$-coloring for $c\in\{1,2\}$. 
The situation is slightly different for $c\in\{3,4\}$, since a $5$-cycle is a graph with maximum average degree exactly 2 and requires five colors in a proper conflict-free coloring. 
Nonetheless, every graph $G$ with $\mad(G)<2$ is a tree, which has a proper conflict-free 3-coloring, demonstrated by picking an arbitrary vertex $v$ and coloring every vertex of the tree by its distance modulo 3 from $v$. 

Let $K_c$ be the complete graph on $c$ vertices, and let $K^*_{c}$ be the graph obtained by subdividing every edge of $K_c$ exactly once. 
For $c\geq 5$, we prove the following sharp threshold on the maximum average degree that guarantees a proper conflict-free $c$-coloring:

\begin{theorem}\label{thm:mad:c6}
For $c\geq 5$,  if $G$ is a graph with $\mad(G)\leq\mad(K^*_{c+1})=\frac{4c}{c+2}$, then $\pcf(G)\leq c$, unless $G$ contains $K^*_{c+1}$ as a subgraph.
\end{theorem}

Similarly, for $c=4$, it turns out the 5-cycle is the only tightness example, as forbidding the 5-cycle allows us to bypass maximum average degree 2: 

\begin{theorem}\label{thm:mad:c4}
If $G$ is a graph with $\mad(G)<\frac{12}{5}$ and no induced $5$-cycle, then $\pcf(G)\leq 4$.
\end{theorem}

There is no hope for an analogous result for $c=3$, as every cycle on $3k\pm1$ vertices does not have a proper conflict-free 3-coloring.

Theorem~\ref{thm:mad:c6} implies that every planar graph with girth at least 5 has a PCF 10-coloring.
We make the following improvement: 

\begin{theorem}\label{thm:g5c7}
If $G$ is a planar graph with girth at least $5$, then $\pcf(G)\leq 7$.
\end{theorem}

We prove all results via the discharging technique.
In Section~2, we collect all reducible configurations utilized in the proofs of Theorem~\ref{thm:mad:c6} and Theorem~\ref{thm:mad:c4}.
In Section~3, the discharging procedures for Theorem~\ref{thm:mad:c6} and Theorem~\ref{thm:mad:c4} are laid out. 
We had to split Section~3 into three subsections, depending on the number of colors. 
The proof of Theorem~\ref{thm:g5c7} is in Section~\ref{sec:g5c7}.

We end the introduction with some frequently used notation. 
A $d$-vertex (resp. $d^-$- and $d^+$-vertex) is a vertex of degree exactly $d$ (resp. at most $d$ and at least $d$). 
A $d$-neighbor (resp. $d^-$- and $d^+$-neighbor) of a vertex $v$ is a $d$-vertex (resp. $d^-$- and $d^+$-vertex) that is a neighbor of $v$.
For a vertex $v$, we denote by $d_G(v)$ and $N_G(v)$ the degree and the neighborhood, respectively, of $v$. 
We also define $N_{G,d}(v)$ (resp. $N_{G,d^-}(v)$ and $N_{G,d^+}(v)$) to be the set of all $d$-neighbors (resp. $d^-$- and $d^+$-neighbors) of $v$, and 
$n_{G,d}(v)=|N_{G,d}(v)|$  (resp. $n_{G,d^-}(v)=|N_{G,d^-}(v)|$ and $n_{G,d^+}(v)=|N_{G,d^+}(v)|$).
If there is no confusion, then we drop $G$ in our notation, as in  $d(v)$, $N(v)$, $N_d(v)$, $N_{d^+}(v)$, $n_d(v)$, $n_{d^+}(v)$ and so on. 
A {\it $k$-thread} means a path on  $k$ $2$-vertices.

\section{Graphs with bounded maximum average degree: Reducible Configurations}

Let $G_c$ be a graph with no PCF $c$-coloring where every proper subgraph of $G_c$ has a PCF $c$-coloring. 
In the proof of each reducible configuration, we will define a non-empty set $S$ of vertices and obtain a PCF $c$-coloring $\vp$ of $G_c-S$, which is guaranteed to exist by the minimality of $G_c$. 
We then extend $\vp$ to $G_c$, which is a contradiction, to conclude that the desired configuration does not appear in $G_c$.
For a vertex $v$, let $\uc(v)$ denote a color that appears only once on the neighborhood of $v$, if it exists.  We will consider $\uc(v)$ to be automatically updated as the coloring of the graph changes. 
In addition, for technical reasons, we allow a PCF $c$-coloring of the null graph, which is a graph whose vertex set is empty.

\begin{lemma} \label{lem:4rc-basic}
For $c\ge 4$, the following do not appear in $G_c$:
\begin{enumerate}[(i)]
    \item A $1$-vertex.
    \item A $3$-vertex adjacent to a $2$-thread.
   \item A $4$-thread not on an induced $5$-cycle. 
\end{enumerate}
\end{lemma}

\begin{proof} 
(i) Suppose to the contrary that $G_c$ has a $1$-vertex $v$, and let $u$ be the neighbor of $v$. 
Let $S=\{v\}$, and obtain a PCF $c$-coloring $\vp$ of $G_c-S$. 
Use a color not in $\{\vp(u), \uc(u)\}$ on $v$ to extend $\vp$ to $G_c$.
 
(ii)
Suppose to the contrary that $G_c$ has a path $u_1v_1v_2u_2$ where $v_1v_2$ is a $2$-thread and $u_1$ is a $3$-vertex. 
Let $S=\{v_1, v_2\}$, and obtain a PCF $c$-coloring $\vp$ of $G_c-S$.
Note that $\uc(u_1)$ is defined at this point. 
Use a color not in $\{\vp(u_2), \uc(u_2), \vp(u_1)\}$ on  $v_2$, then  use a color not in $\{\vp(v_2), \vp(u_2), \vp(u_1)\}$ on $v_1$ to extend $\vp$ to  $G_c$.

(iii) Suppose to the contrary that $G_c$ has a path $u_1v_1v_2v_3v_4u_2$ where $v_1v_2v_3v_4$ is a $4$-thread and $u_1\neq u_2$.
Let $S=\{v_1, v_2, v_3, v_4\}$, and obtain a PCF  $c$-coloring $\vp$ of $G_c-S$. 

If $\{\vp(u_1), \uc(u_1)\}$ and $\{\vp(u_2), \uc(u_2)\}$ are not disjoint,  then there exists a color that can be used on both $v_1$ and $v_4$. 
Now,  greedily color $v_2$ and $v_3$ to extend $\vp$ to  $G_c$. 

If $\{\vp(u_1), \uc(u_1)\}$ and $\{\vp(u_2), \uc(u_2)\}$ are disjoint,
then color $v_1$ with $\vp(u_2)$ and $v_4$ with $\vp(u_1)$.
Now,  greedily color $v_2$ and $v_3$ to extend $\vp$ to  $G_c$. 
\end{proof}

\begin{lemma}\label{lem:thread}
For $c \ge 4$ and  $d\in \{4, 5\}$, if a $d$-vertex 
$v$ of $G_c$ is adjacent to a $3$-thread, then 
$v$ is adjacent to at most $3d + n_{3^+}(v)-11$  $2$-threads. 
\end{lemma}
\begin{proof}
Let $uxy$ be a $3$-thread where $u$ is adjacent to $v$.
Let $T'$ be the set of $2$-neighbors of $v$ that are not on a $2$-thread, $T''$ be the set of $2$-vertices on a $2$-thread adjacent to $v$, and $T=T'\cup T''$. 
Note that the number of $2$-threads adjacent to $v$ is exactly $n_2(v)-|T'|$.
Suppose to the contrary that $v$ is adjacent to at least $3d+n_{3^+}(v)-10$ $2$-threads, so
 $2n_{3^+}(v)+|T'|\le 10-2d$.

Let $S=T\cup\{v,y\}$,
and obtain a PCF $c$-coloring $\vp$ of $G_c-S$. 
Let $C = \vp(N_{3^+}(v) \cup (N(T')\setminus\{v\})) \cup  \uc(N_{3^+}(v))$, so $|C|\le 10-2d$. 
Namely, if $d=5$ then $|C|=0$, and if $d=4$ then $|C|\le 2$. 

Let $z$ be the neighbor of $y$ that is not $x$. 
In each case, color the vertex $v$ as the following:
\begin{enumerate}
    \item If $d=5$, then use the color $\vp(z)$ on $v$; if $\vp(z)$ is not defined (which is when $z=v$), then use an arbitrary color on $v$.
    \item If $d=4$ and $\uc(z)$ is defined, then use the color $\uc(z)$ on $x$, and use a color not in $C \cup \{\vp(x)\}$ on $v$.
    \item If $d=4$ and $\uc(z)$ is not defined (which is when no neighbors of $z$ are colored yet), then use a color not in $C \cup \{\vp(z)\}$ on both $v$ and $y$.
\end{enumerate}
Color the vertices in $T\setminus\{u,x,y\}$ as follows. 
Let $v_1\ldots v_t$ be a $t$-thread in $T$ where $v_1$ is adjacent to $v$ and let $v_{t+1}$ be the neighbor of $v_t$ that is not $v_{t-1}$. 
Color $v_t, \ldots, v_1$ in this order, and use a color not in $\{\vp(v_{i+1}), \uc(v_{i+1}), \vp(v)\}$ on $v_i$.

In case of (1), if $\uc(v)$ is not  defined, then at most two colors $\alpha$ and $\beta$ appear on the neighbors of $v$ that are not $u$.
Otherwise, let $\alpha = \uc(v)$.
Now, use a color not in $\{\vp(v),\alpha,\beta\}$ on $u$, then use a color not in $\{\vp(z), \uc(z), \vp(u)\}$ on $y$, then use a color not in $\{\vp(v), \vp(u), \vp(y)\}$ on $x$ to extend $\vp$ to   $G_c$.

In case of (2), if $\uc(v)$ is not defined, then only one color $\alpha$ appears on the neighbors of $v$ that are not $u$.
Otherwise, let $\alpha = \uc(v)$.
Now, use a color not in $\{\vp(v),\alpha, \vp(x)\}$ on $u$, then use a color not in $\{\vp(z), \uc(z), \vp(u)\}$ on $y$ to extend $\vp$ to $G_c$.

In case of (3), define $\alpha$ as in the case of (2).
Now, use a color not in $\{\vp(v), \alpha\}$ on $u$, then use a color not in $\{\vp(v),\vp(u),\vp(z)\}$ on $x$ to extend $\vp$ to  $G_c$.
\end{proof}

\begin{lemma}\label{lem:combine} 
For $c\ge 5$, if $v_1$ and $v_2$ are adjacent $3^+$-vertices in $G_c$ such that $2d(v_i) - n_2(v_i) -2 \le c-i$ for each $i \in \{1,2\}$, then $v_1$ is a $4^+$-vertex, and either $v_1$ or $v_2$ has no $2$-neighbors.
\end{lemma}

\begin{proof}
Consider a PCF $c$-coloring $\vp$ of $G_c-N_2(v_1) - N_2(v_2)-\{v_1,v_2\}$. 
Let $C_i=\vp(N_{3^+}(v_i) \setminus\{v_{3-i}\}) \cup \uc((N_{3^+}(v_i) \setminus\{v_{3-i}\}) \cup N_2(v_i))$ for each $i \in \{1,2\}$.
Note that $|C_i| \le 2(d(v_i)-n_2(v_i)-1) + n_2(v_i) = 2d(v_i) - n_2(v_i) -2 \le c-i$ for each $i \in \{1,2\}$.

\begin{claim}\label{claim:sub}
Every extension of $\vp$ such that $\vp(v_1)$, $\vp(v_2)$, $\uc(v_1)$, and $\uc(v_2)$ are defined can be further extended to a PCF $c$-coloring of all of $G$.
\end{claim}
\begin{proof} 
If $\vp(v_1)$, $\vp(v_2)$, $\uc(v_1)$, and $\uc(v_2)$ are defined, then extend $\vp$ to  $G$ by coloring the vertices in $N_2(v_1)\cup N_2(v_2)$ one by one: use a color not in $\vp(N(v))\cup\uc(N(v))$ on $v$ for each $v \in N_2(v_1)\cup N_2(v_2)$.
Note that $|\vp(N(v))\cup\uc(N(v))|\leq 4\leq c-1$ for each $v\in N_2(v_1)\cup N_2(v_2)$. 
\end{proof}

Suppose to the contrary that for each $i\in\{1,2\}$, $v_i$ has a $2$-neighbor $u_i$, whose neighbor other than $v_i$ is $x_i$.
Use a color not in $C_1$ on $v_1$, then use a color not in $C_2 \cup \{\vp(v_1)\}$ on $v_2$.
If $\uc(v_i)$ is not defined, then use a color not in $\vp(N_{3^+}(v_i))\cup\{\vp(v_i),\vp(x_i),\uc(x_i)\}$ on $u_i$.
Note that  
\[ |\vp(N_{3^+}(v_i))|+3\le \left\lfloor\frac{d(v_i)-n_2(v_i)}{2}\right\rfloor+3 \le \left\lfloor\frac{c-n_2(v_i)+1}{4}\right\rfloor +3 \le  c-1,\]
where the first inequality holds since $\uc(v_i)$ is not defined, the second inequality holds since $2d(v_i) - n_2(v_i) -2 \le c-1$, and the third inequality holds since $c \ge 5$. 
Hence, we may assume both $\uc(v_1)$ and $\uc(v_2)$ are defined, so by Claim~\ref{claim:sub}, $\vp$ can be extended to  $G$.

Now suppose to the contrary that $v_1$ is a $3$-vertex with neighbors $u_1,u_2,v_2$. 
If $\uc(v_2)$ is not defined, then color $v_1$ with a color not in $C_1^*=C_1 \cup \vp(N_{3^+}(v_2) \setminus \{v_1\})$. Note that  \begin{eqnarray*}
|C_1^*| \le  |C_1| +|\vp(N_{3^+}(v_2) \setminus \{v_1\})| \le |C_1|+ \left\lfloor \frac{d(v_2)-n_2(v_2)-1}{2} \right\rfloor  &\le&  |C_1| + \left\lfloor \frac{c-2-n_2(v_2)}{4} \right\rfloor \\&\le& |C_1|+c-5,\end{eqnarray*} 
where the second inequality holds since $\uc(v_2)$ is not defined, the third inequality holds since $2d(v_2)-n_2(v_2)-2 \le c-2$, and the fourth inequality holds since $c \ge 5$.
Since $|C_1|\leq 4$, we can color $v_1$ since $|C^*_1|\leq c-1$. 
Hence, we may assume $\uc(v_2)$ is defined.

If $\uc(v_1)$ is not defined, then $\vp(u_1)=\vp(u_2)$, so $|C_1|\le 3$, which further implies $|C^*_1|\leq c-2$. 
Therefore, there are two colors not in $C_1^*$ that can be used on $v_1$.  
Now, color $v_2$ with a color not in $C_2\cup \{ \vp(u_1)\}$, then recolor $v_1$ with a color not in $C_1^*\cup\{ \vp(v_2)\}$. 
So we may assume $\uc(v_1)$ is defined. 

Since $\vp(v_1)$, $\vp(v_2)$, $\uc(v_1)$ and $\uc(v_2)$ are defined, by Claim~\ref{claim:sub} $\vp$ can be extended to $G$.
\end{proof}
 
\begin{lemma}\label{lem:2neighbor-new}
Let $v$ be a vertex with a $2$-neighbor.
If $c\ge 5$, then $2d(v)\geq n_2(v)+c$.
Moreover, if $c\ge 7$, then $2d(v)\geq n_2(v)+n_3(v)+c$.
\end{lemma}
\begin{proof} 
Suppose to the contrary that $2d(v) \le n_2(v)+c-1$ when $c \ge 5$, and $2d(v) \le n_2(v) + n_3(v) + c -1$ when $c \ge 7$.
 For $i \in \{2,3\}$, let $N'_i(v)$ be the set of vertices consisting of a neighbor $u$ for each vertex in $N_i(v)$ where $u \neq v$.

When $c \ge 5$, let $\vp$ be a PCF $c$-coloring of $G_c - (N_2(v) \cup \{v\})$, and color $v$ with a color not in $C= \vp(N'_2(v) \cup N_{3^+}(v)) \cup \uc(N_{3^+}(v))$.
Note that \[|C| \le n_2(v) + 2n_{3^+}(v) \le n_2(v) + 2(d(v)-n_2(v)) \le c-1. \]
When $c \ge 7$, let $\vp$ be a PCF $c$-coloring of $G_c-(N_2(v) \cup N_3(v) \cup \{v\})$, and color $v$ with a color not in $C'= \vp(N'_2(v) \cup N'_3(v) \cup N_{4^+}(v)) \cup \uc(N_{4^+}(v))$. 
Note that \[ |C'| \le n_2(v)+n_3(v)+2n_{4^+}(v)\le n_2(v)+n_3(v) + 2(d(v)-
n_2(v)-n_3(v))   \le c-1.\]
In both cases, we will extend $\vp$ to  $G$ by coloring the vertices in $N_2(v)$ (resp. $N_2(v)\cup N_3(v)$) when $c \ge 5$ (resp. $c \ge 7$). 
For a particular vertex $u_0 \in N_2(v)$, color $u_0$ with a color not in $C_1 = \vp(N_{3^+}(v)) \cup \{\vp(v), \vp(u_0'), \uc(u_0')\}$ (resp. $C_{1}'=\vp(N_{4^+}(v)) \cup \{\vp(v), \vp(u_0'), \uc(u_0')\}$) where $u_0'$ is the neighbor of $u_0$ that is not $v$, so $\uc(v)=\vp(u_0)$. 
Note that
\[|C_1|\le n_{3^+}(v)+3 \le  
\left \lfloor\frac{c-1-n_2(v)}{2} \right\rfloor +3   \le  c-1,\]
and
\[|C'_1|\le n_{4^+}(v)+3 \le  
\left \lfloor\frac{c-1-n_2(v)-n_3(v)}{2} \right\rfloor +3   \le  c-1,\] as $c \ge 5$ and $n_2(v)\ge 1$.
For a vertex $u\in N_2(v)\setminus \{u_0\}$, use a color not in $\{\vp(u'), \uc(u'), \vp(v), \vp(u_0)\}$ on $u$, where $u'$ is the neighbor of $u$ that is not $v$. 
If $c\ge 7$, then for each $w \in N_3(v)$, use a color not in $\vp(N(w))\cup\uc(N(w))$ on $w$.
Note that $|\vp(N(w))\cup\uc(N(w))|\leq 6\le c-1$. 
\end{proof}

\section{Graphs with bounded maximum average degree: Discharging Rules}

In this section, we prove both Theorem~\ref{thm:mad:c6} and Theorem~\ref{thm:mad:c4} by the discharging method.
For each theorem, let $G$ be a counterexample with the minimum number of vertices. 
For each vertex $v \in V(G)$, let its initial charge $\mu(v)$ be the degree of $v$.
Via a set of discharging rules the charge will be redistributed, to obtain final charge $\mu^*(v)$ at each vertex $v$. 
The total charge sum will be preserved during the discharging process.
Either the final charge sum will be different from the initial charge sum or we conclude that $G$ is a graph that actually has a PCF $c$-coloring, to conclude that a counterexample could not have existed in the first place. 

We introduce subsections as the discharging rules vary slightly depending on the number of colors. 
The discharging process for Theorem~\ref{thm:mad:c4} is in  subsection~\ref{sub:c4}.
The discharging process for Theorem~\ref{thm:mad:c6} is split into two subsections, depending on the number of colors. 
The number of colors is exactly five and at least six in subsection~\ref{sub:c5} and subsection~\ref{sub:c6}, respectively.

\subsection{Four colors}\label{sub:c4}

Let $G$ be a counterexample to Theorem~\ref{thm:mad:c4} with the minimum number of vertices.
In other words, $G$ is a graph with $\mad(G) < \frac{12}{5}$, $G$ has no induced $5$-cycle, and $G$ has no PCF $4$-coloring, while every proper subgraph of $G$ has a PCF $4$-coloring.

A $2$-vertex $u$ on a thread adjacent to a $3^+$-vertex $v$ is {\it close} to $v$, and $u$ is a {\it close neighbor} of $v$.
There is one discharging rule:

\begin{enumerate}[(R1)]
    \item\label{madrule:1} 
    Each $3^+$-vertex sends charge $\frac{1}{5}$ to each of its close neighbors. 
\end{enumerate}

Consider a vertex $v$.  
Note that $v$ is not adjacent to a $4$-thread by Lemma~\ref{lem:4rc-basic}~(iii) since $G$ has no induced $5$-cycle.

Note that $G$ has no $1$-vertex by Lemma~\ref{lem:4rc-basic}~(i).
If $v$ is a $2$-vertex, then it receives charge $\frac{1}{5}$ from each vertex $u$ where $v$ is a close neighbor of $u$ by~(R\ref{madrule:1}).
Hence, $\mu^*(v) = 2 + 2 \cdot \frac{1}{5} = \frac{12}{5}$.
If $v$ is a $3$-vertex, then by Lemma~\ref{lem:4rc-basic} (ii), $v$ has at most three close neighbors. 
Hence, $\mu^*(v) \geq 3 - 3 \cdot \frac{1}{5} = \frac{12}{5}$.

Suppose $v$ is a $d$-vertex where $d\in\{4,5\}$. 
If $v$ is not adjacent to a $3$-thread, then $v$ has at most $2d$ close neighbors.
Hence, $\mu^*(v) \geq d - 2d \cdot \frac{1}{5} \geq \frac{12}{5}$.
If $v$ is adjacent to a $3$-thread, then by Lemma~\ref{lem:thread}, 
$d \ge n_{3^+}(v)+ t\ge 11-3d + 2t$,
where $t$ is the number of $2$-threads adjacent to $v$.
Thus, $t \le \left\lfloor \frac{4d-11}{2} \right\rfloor$, so $v$ has at most $d+2\left\lfloor \frac{4d-11}{2} \right\rfloor$ close neighbors. 

Hence, $\mu^*(v) \geq d - \left(d+2\left\lfloor \frac{4d-11}{2} \right\rfloor\right)\cdot\frac{1}{5} = \frac{12}{5}$.
If $v$ is a $6^+$-vertex, then $v$ has at most $3d(v)$ close neighbors.
Hence, $\mu^*(v) \ge d(v)-3d(v) \cdot \frac{1}{5}  \ge \frac{12}{5}$ since $d(v) \ge 6$.

Thus, $\mu^*(v) \geq \frac{12}{5}$ for every vertex $v$ in $G$, so the final charge sum is at least $\frac{12}{5}|V(G)|$.
This is a contradiction since the initial charge sum is $\sum_{v\in V(G)}d(v)=2|E(G)|<\frac{12}{5}|V(G)|$. 

\subsection{Five colors}\label{sub:c5}

For $c=5$, let $G$ be a counterexample to Theorem~\ref{thm:mad:c6} with the minimum number of vertices.
In other words, $G$ is a graph with $\mad(G) \leq \frac{20}{7}$, $G$ does not have $K_6^*$ as a subgraph, and $G$ has no PCF $5$-coloring, while every proper subgraph of $G$ has a PCF $5$-coloring.

The discharging rules are the following:
\begin{enumerate}[(R1)]
    \item\label{rule:c=5_1} 
    Each $3^+$-vertex sends charge $\frac{3}{7}$ to each of its $2$-neighbors.
    \item\label{rule:c=5_2} 
    Each $4^+$-vertex sends charge $\frac{1}{7}$ to each of its $3$-neighbors $u$ with $n_2(u) \ge 1$, and to each of its $4$-neighbors $w$ with $n_2(w) \ge 3$.
\end{enumerate}

Consider a vertex $v$.
We will first show that $\mu^*(v) \geq \frac{20}{7}$.
Note that $G$ has no $1$-vertex by Lemma~\ref{lem:4rc-basic}~(i).
If $v$ is a $2$-vertex, then every neighbor of $v$ is a $3^+$-vertex by Lemma~\ref{lem:2neighbor-new}. 
By (R\ref{rule:c=5_1}), $\mu^*(v) = 2+ 2 \cdot \frac{3}{7} = \frac{20}{7}$.

If $v$ is a $3$-vertex, then $v$ sends charge to only $2$-neighbors.
If $v$ has no $2$-neighbors, then $v$ is not involved in the discharging process, so $\mu^*(v) = 3 > \frac{20}{7}$. 
If $v$ has a $2$-neighbor, then $n_3(v)=0$ by Lemma~\ref{lem:combine}, and Lemma~\ref{lem:2neighbor-new} further implies $v$ has exactly one $2$-neighbor.
Since $v$ also has two $4^+$-neighbors, by the rules, $\mu^*(v) = 3 - \frac{3}{7} + 2 \cdot \frac{1}{7} = \frac{20}{7}$.

If $v$ is a $4$-vertex, then $v$ has at most three $2$-neighbors by Lemma~\ref{lem:2neighbor-new}.
If $v$ has at most one $2$-neighbor, then by the rules, $\mu^*(v) \ge 4 - \frac{3}{7} - 3 \cdot \frac{1}{7} > \frac{20}{7}$.
If $v$ has either two or three $2$-neighbors, then Lemma~\ref{lem:combine} implies $v$ has neither a $3$-neighbor $u$ with $n_2(u)\ge 1$ nor a $4$-neighbor $w$ with $n_2(w) \ge 3$.
Hence, by the rules, $\mu^*(v) \ge 4 - \max\{2 \cdot \frac{3}{7}, 3 \cdot \frac{3}{7} - \frac{1}{7}\}\geq \frac{20}{7}$.

If $v$ is a $5^+$-vertex, then by the rules, $\mu^*(v) \ge d(v) - d(v) \cdot \frac{3}{7} = \frac{4}{7}d(v) \ge \frac{20}{7}$.

Thus, $\mu^*(v) \ge \frac{20}{7}$ for every vertex $v$ in $G$.
Since the initial charge sum is at most $\frac{20}{7}|V(G)|$, we obtain that $\mu^*(v) = \frac{20}{7}$ for every vertex $v$ in $G$.  
By observing the final charge of each vertex of $G$, each vertex of $G$ is either a 2-vertex or a $d$-vertex $v$ with exactly $2d-5$ $2$-neighbors for some $d\in\{3,4,5\}$.

Since a $5$-vertex has only $2$-neighbors,
a $3$-vertex or a $4$-vertex must have a $3$-neighbor or a $4$-neighbor, which is impossible by Lemma~\ref{lem:combine}.
Thus, every $3^+$-vertex of $G$ is a $5$-vertex with exactly five $2$-neighbors. 
Therefore, $G$ is a graph obtained from a $5$-regular graph $G_0$ by subdividing every edge exactly once.
Since $G$ does not have $K_6^*$ as a subgraph, $G_0$ does not have $K_6$ as a subgraph, and thus $G_0$ is properly $5$-colorable by Brook's Theorem \cite{1941Brooks}. 
Let $\vp$ be a proper $5$-coloring of $G_0$, and consider $\vp$ as a coloring of the $5$-vertices of $G$.
Color the $2$-vertices of $G$ one by one as follows.
For each $2$-vertex $v$, use a color not in $\vp(N(v)) \cup \uc(N(v))$; this is possible since $|\vp(N(v)) \cup \uc(N(v))|\leq 4$.
Now, $\vp$ is a PCF $5$-coloring of $G$, which is a contradiction.

\subsection{At least six colors}\label{sub:c6}
 
For $c\geq6$, let $G$ be a counterexample to  Theorem~\ref{thm:mad:c6} with the minimum number of vertices.
In other words, $G$ is a graph with $\mad(G)\leq \frac{4c}{c+2}$, $G$ does not have $K^*_{c+1}$ as a subgraph, and $G$ has no PCF $c$-coloring, while every proper subgraph of $G$ has a PCF $c$-coloring.

The discharging rules are the following:
\begin{enumerate}[(R1)]
    \item\label{madrule:1} 
    Each $4^+$-vertex sends charge $\frac{c-2}{c+2}$ to each of its $2$-neighbors.
    \item\label{madrule:2} 
    Each $4^+$-vertex sends charge $\frac{c-6}{3(c+2)}$ to each of its $3$-neighbors.
\end{enumerate} 
Consider a vertex $v$.  
We will first show that $\mu^*(v) \geq \frac{4c}{c+2}$.
Note that $G$ has neither a $1$-vertex nor a $2$-thread by Lemma~\ref{lem:4rc-basic}~(i) and Lemma~\ref{lem:2neighbor-new}, respectively.
In addition, a $3$-vertex $v$ has no $3^-$-neighbor by Lemma~\ref{lem:combine} and Lemma~\ref{lem:2neighbor-new}.
Thus, a $3^-$-vertex has no $3^-$-neighbor.

If $v$ is a $2$-vertex, then it receives charge $\frac{c-2}{c+2}$ from each neighbor by~(R\ref{madrule:1}) since $v$ has only $4^+$-neighbors.
Hence, $\mu^*(v) = 2 + 2 \cdot \frac{c-2}{c+2} = \frac{4c}{c+2}$.
If $v$ is a $3$-vertex, then it receives charge $\frac{c-6}{3(c+2)}$ from each neighbor by (R\ref{madrule:2}) since $v$ has only $4^+$-neighbors.
Hence, $\mu^*(v)\ge 3+3\cdot\frac{c-6}{3(c+2)}=\frac{4c}{c+2}$.

In the following, suppose that $v$ is a $4^+$-vertex. 
If $c=6$, then $v$ sends charge only to its $2$-neighbors and none to its $3$-neighbors. 
Lemma~\ref{lem:2neighbor-new} implies $n_2(v)\le 2d(v) -6$, so $\mu^*(v)\ge d(v)-(2d(v)-6)\cdot\frac{c-2}{c+2}\geq \frac{4c}{c+2}$.

Now suppose $c\geq 7$. 
The discharging rules imply $\mu^*(v)\ge d(v)-n_2(v)\cdot \frac{c-2}{c+2}-n_3(v)\cdot\frac{c-6}{3(c+2)}$.
Since  $d(v)=n_{4^+}(v)+n_3(v)+n_2(v)$, the following holds:
\begin{eqnarray}\label{eq1}&&
\mu^*(v) \geq   n_{4^+}(v)+ \frac{2c+12}{3(c+2)}\cdot n_3(v)  +\frac{4}{c+2}\cdot n_2(v).  
\end{eqnarray}
Since $1\ge \frac{2c+12}{3(c+2)}\ge \frac{4}{c+2}$, 
\eqref{eq1} implies $\mu^*(v)\ge d(v)\cdot \frac{4}{c+2}$.
If either $d(v)\ge c$ or $n_2(v)+n_3(v)=0$, then $\mu^*(v)\ge \min\{4,\frac{4c}{c+2}\}=\frac{4c}{c+2}$. 
Thus we may assume that $d(v)\leq c-1$ and $n_2(v)+n_3(v)\ge 1$.

First, suppose that $n_2(v)=0$, so $n_3(v)>0$. 
If $d(v)\in \{4,5\}$ and $c\ge 10$, then 
applying Lemma~\ref{lem:combine} to $v$ and a $3$-neighbor of $v$ results in the following contradiction:  $8\geq 2d(v)-2\ge c-1 \ge 9$. 
Thus  $d(v)\ge 6$ or $c\in\{7,8,9\}$,  so 
\eqref{eq1} implies  $\mu^*(v)\ge d(v)\cdot \frac{2c+12}{3(c+2)}>\frac{4c}{c+2}$.

Now suppose that $n_2(v) > 0$.
Since $c \ge 7$, Lemma~\ref{lem:2neighbor-new} implies $n_2(v) + n_3(v) \le 2d(v) -c$, so \eqref{eq1} implies
\[
\mu^*(v)\ge 
(c-d(v)) +\frac{4}{c+2}\cdot (2d(v)-c)  \ge \frac{c^2-2c - (c-6)d(v)}{c+2} \ge  \frac{5c-6}{c+2} >\frac{4c}{c+2},
\]
where the penultimate inequality holds since $c\ge 7$ and $d(v)\le c-1$.

Thus $\mu^*(v) \ge \frac{4c}{c+2}$ for every vertex $v$ of  $G$.
Since the initial charge sum is at most $\frac{4c}{c+2}|V(G)|$, we obtain that $\mu^*(v) = \frac{4c}{c+2}$ for every vertex $v$ in $G$.  

If $c = 6$, by observing the final charge of each vertex of $G$, 
each $3^+$-vertex of $G$ is actually a $d$-vertex with exactly 
$2d-6$ $2$-neighbors for $d \in \{3,4,5,6\}$.
Since a $6$-vertex has only $2$-neighbors, for $d \in \{3,4,5\}$, a $d$-vertex must have a $d'$-neighbor for some $d' \in \{3,4,5\}$, which is impossible by Lemma~\ref{lem:combine}.
Thus, every $3^+$-vertex of $G$ is a $6$-vertex with exactly six $2$-neighbors.

If $c \ge 7$, by observing the final charge of each vertex of $G$, each vertex of $G$ is either a $c$-vertex with exactly $c$ $2$-neighbors or a $3^-$-vertex.
Since a $3^-$-vertex has no $3^-$-neighbor  
there is no vertex in $G$ that is a neighbor of a $3$-vertex. 
Thus, every $3^+$-vertex of $G$ is a $c$-vertex with exactly $c$ $2$-neighbors.

In each case, $G$ is obtained from a $c$-regular graph $G_0$ by subdividing every edge exactly once.
Since $G$ does not have $K_{c+1}^*$ as a subgraph, $G_0$ does not have $K_{c+1}$ as a subgraph, and thus $G_0$ is properly $c$-colorable by Brooks' Theorem~\cite{1941Brooks}.
Let $\vp$ be a proper $c$-coloring of $G_0$, and consider $\vp$ as a coloring of the $c$-vertices of $G$. 
Color the 2-vertices of $G$ one by one as follows. 
For each $2$-vertex $v$, use a color not in $\vp(N(v))\cup\uc(N(v))$; this is possible since $|\vp(N(v) \cup \uc(N(v))| \le 4$.
Now, $\vp$ is a PCF $c$-coloring of $G$, which is a contradiction.

\section{Planar graphs with girth at least $5$ are PCF $7$-colorable}
\label{sec:g5c7}

Let $G$ be a counterexample to Theorem~\ref{thm:g5c7} with the minimum number of vertices, and fix a plane embedding of $G$.
So, $G$ is a plane graph with girth at least $5$ that has no PCF $7$-coloring, while every proper subgraph of $G$ has a PCF $7$-coloring.

\begin{figure}[h!]
    \centering
    \includegraphics[height=3.5cm]{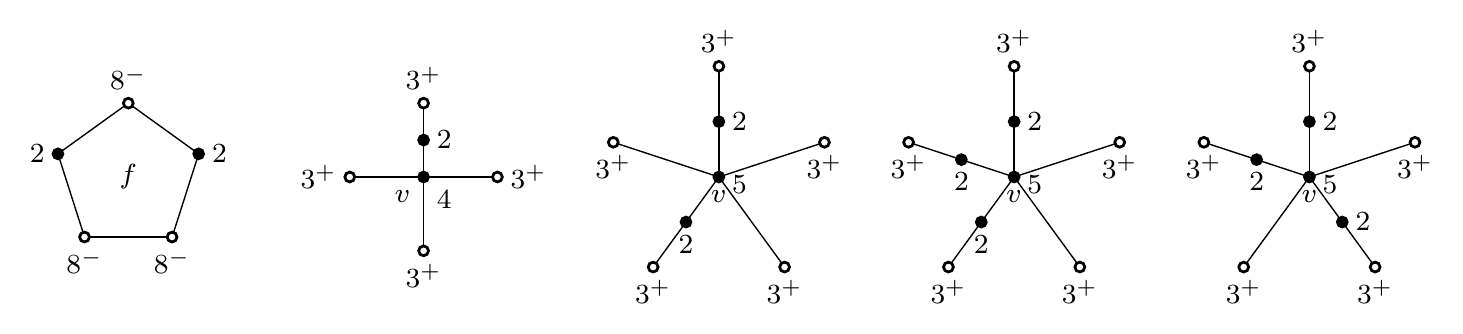}
    \caption{Illustrations for a terrible face $f$ and a bad vertex $v$.}
    \label{fig:tfandbv}
\end{figure}

A {\it $d$-face} is a face of length $d$. 
For an edge $uv$ on the boundary of a face $f$, the vertex $u$ is a {\it boundary neighbor} of $v$. 
A $4^+$-vertex $v$ is \textit{giving} on $f$ if $v$ has a boundary $3^+$-neighbor on $f$. 
For a $4^+$-vertex $v$, let $F^*_{5}(v)$ denote the set of $5$-faces $f$ incident with $v$ such that $v$ is giving on $f$. 
If a $5$-face $f$ is incident with two $2$-vertices and three $8^-$-vertices, then $f$ is a \textit{terrible} face.
Note that the $2$-vertices on a terrible face are not adjacent to each other by  Lemma~\ref{lem:2neighbor-new}.
A vertex $v$ is \textit{bad} if $v$ is one of the following:
\begin{itemize}
    \item a $4$-vertex with exactly one $2$-neighbor;
    \item a $5$-vertex with exactly two $2$-neighbors and $|F^*_5(v)|=5$;
    \item a $5$-vertex with exactly three $2$-neighbors.
\end{itemize}
A $4^+$-vertex is \textit{good} if it is not bad.  
See Figure~\ref{fig:tfandbv} for illustrations of a terrible face and a bad vertex.

For a vertex $v$, let $t(v)$ denote the number of boundary $3^+$-neighbors of $v$ on a terrible face incident with $v$.  
Recall that there are no adjacent $3^-$-vertices by Lemma~\ref{lem:2neighbor-new}.

\begin{lemma}\label{lem:t(v)}
If $v$ is a bad vertex, then $2d(v)-n_2(v)-n_3(v)-t(v) \ge 7$.
\end{lemma}
\begin{proof}
Suppose to the contrary that $2d(v)-n_2(v)-n_3(v)-t(v) \le 6$.
Let $X$ be the set of $3^+$-neighbors of $v$ on terrible faces incident with $v$, so  $|X|=t(v)$.
Note that $X \subseteq N_{4^+}(v)$ since there are no adjacent $3^-$-vertices. 
Let $S= N_2(v) \cup N_3(v) \cup \{v\}$, and obtain a PCF $7$-coloring $\vp$ of $G-S$.
For $i \in \{2,3\}$, let $N'_i(v)$ be the set of vertices consisting of a neighbor $u$ for each vertex in $N_i(v)$ where $u \neq v$.
Use a color not in $\vp(N'_2(v) \cup N'_3(v) \cup N_{4^+}(v)) \cup \uc(N_{4^+}(v) \setminus X)$ on $v$, which is possible since $n_2(v) + n_3(v) + t(v) + 2(d(v)-n_2(v)-n_3(v)-t(v)) = 2d(v)-n_2(v)-n_3(v)-t(v) \le 6$.

Let $vu_1u_2u_3u_4v$ be a terrible face incident with $v$, where $u_4 \in X$.
This implies that $u_1$ and $u_3$ are $2$-vertices, and $u_2$ and $u_4$ are $8^-$-vertices.
Erase the color on $u_3$.
Color $u_1$ with a color not in $\vp(N(u_1)) \cup \uc(N(u_1))$.
If $\uc(v)$ is not defined, then all $3^+$-neighbors of $v$ are colored with the same color since $n_{3^+}(v) \le 3$.
If $\uc(u_2)$ is not defined, then there are at most three colors used on the neighbors of $u_2$ since $u_2$ is an $8^-$-vertex. 
Thus there are at most six colors to avoid on $u_1$, so we can always color $u_1$ to guarantee both $\uc(v)$ and $\uc(u_2)$ are defined, if not already. 
Now (re)color $u_3$ with a color not in $\vp(N(u_3)) \cup \uc(N(u_3))$.
If $\uc(u_4)$ is not defined, then there are at most three colors used on the neighbors of $u_4$ since $u_4$ is an $8^-$-vertex. 
Thus there are at most six colors to avoid on $u_3$, so we can always color $u_3$ to guarantee $\uc(u_4)$ is defined, if not already. 

The unique color of a vertex in $X$ is defined by the above argument.
For $u \in N_2(v) \cup N_3(v)$ that is uncolored so far, use a color not in $\vp(N(u)) \cup \uc(N(u))$ on $u$, which is possible since $|\vp(N(u))\cup \uc(N(u))| \le 6$.
\end{proof}

\begin{figure}[h!]
    \centering
    \includegraphics[height=3.5cm]{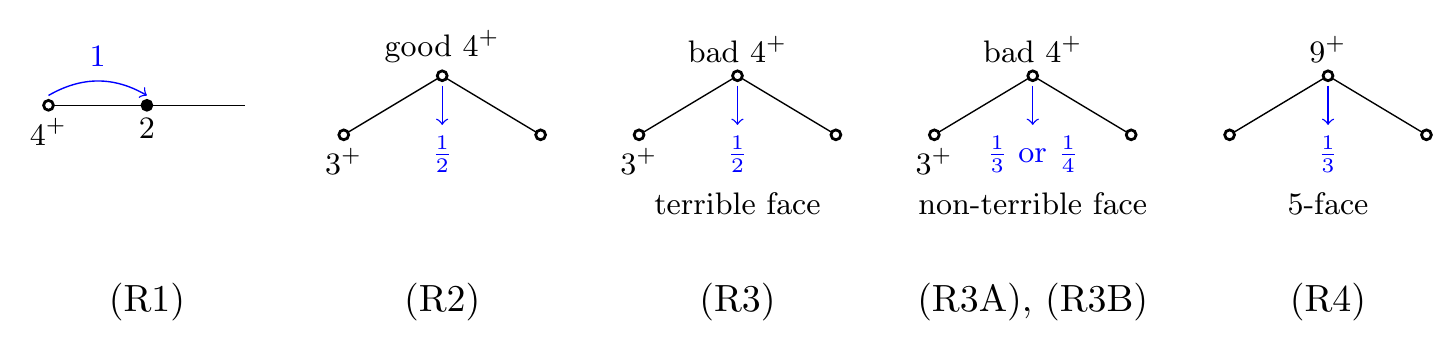}
    \caption{Illustration of the discharging rules.}
    \label{fig:rules:new}
\end{figure}

For a vertex or a face $x$, define the initial charge $\mu(x)$ of $x$ as follows: for each vertex $v$ let $\mu(v) =2d(v) - 6$, and for each face $f$ let $\mu(f) = d(f)-6$. 
Let $\mu^*(x)$ be the final charge of $x$ after distributing the charge by the following discharging rules (see Figure~\ref{fig:rules:new}):

\begin{enumerate}[(R1)]
    \item\label{g5c7rule:2vx:new'} 
    Each $4^+$-vertex sends charge $1$ to each $2$-neighbor.
    \item\label{g5c7rule:good:new'}
    Each good $4^+$-vertex $v$ sends charge $\frac{1}{2}$ to each face in $F^*_{5}(v)$.
    \item\label{g5c7rule:bad:new'}
    Each bad vertex $v$ sends charge $\frac{1}{2}$ to each terrible face in $F^*_5(v)$.
    In addition, 
    \begin{enumerate}[(R3A)]
        \item\label{g5c7rule:bada:new'} if $v$ is a bad $5$-vertex with exactly two $2$-neighbors, then $v$ sends charge $\frac{1}{3}$ to each non-terrible face in $F^*_{5}(v)$;
        \item\label{g5c7rule:badb:new'} otherwise, $v$ sends charge $\frac{1}{4}$ to each non-terrible face in $F^*_{5}(v)$.
    \end{enumerate} 
    \item\label{g5c7rule:9vx:new'}
    Each $9^+$-vertex sends charge $\frac{1}{3}$ to   each  incident $5$-face. 
\end{enumerate}

By Lemma~\ref{lem:2neighbor-new}, a $4$-vertex has at most one $2$-neighbor, and a $5$-vertex has at most three $2$-neighbors. 

\begin{claim}
Each vertex $v$ has non-negative final charge.
\end{claim}
\begin{proof} Note that $G$ has no $1$-vertex by Lemma~\ref{lem:4rc-basic}~(i).
If $v$ is a $2$-vertex, then it has only $4^+$-neighbors by Lemma~\ref{lem:2neighbor-new}, so  $\mu^*(v) = -2 + 2 \cdot 1 =0$ by (R\ref{g5c7rule:2vx:new'}).
If $v$ is a $3$-vertex, then it is not involved in the rules, so $\mu^*(v) =\mu(v) = 0$.

Let $v$ be a $4^+$-vertex. 
Then $v$ sends charge 1 to each $2$-neighbor by (R\ref{g5c7rule:2vx:new'}) and  sends charge at least $\frac{1}{2}$ to each terrible face in $F^*_{5}(v)$ by (R\ref{g5c7rule:good:new'}), (R\ref{g5c7rule:bad:new'}), (R\ref{g5c7rule:9vx:new'}).
Suppose that $v$ is bad. 
If $v$ is incident with a terrible face, then by Lemma~\ref{lem:t(v)}, $v$ is a $5$-vertex with two $2$-neighbors and $t(v)=1$, which further implies that $v$ is incident with at most two terrible faces.
Therefore, if $v$ is a $4$-vertex, then $\mu^*(v) \ge 2 - 1 - 4 \cdot \frac{1}{4}=0$ by (R\ref{g5c7rule:2vx:new'}) and (R\ref{g5c7rule:badb:new'}), and if $v$ is a $5$-vertex, then
$\mu^*(v) \ge 4-\max\{2 \cdot 1 + 2 \cdot \frac{1}{2} + 3 \cdot \frac{1}{3}, 3 \cdot 1 + 4 \cdot\frac{1}{4}\} =0$ by (R\ref{g5c7rule:2vx:new'}) and (R\ref{g5c7rule:bad:new'}).
Note that a $5$-vertex $v$ with three $2$-neighbors is incident with a face not in $F_5^*(v)$

Suppose that $v$ is good. 
If $v$ is a $4$-vertex with no $2$-neighbor or a $5$-vertex with at most one $2$-neighbor, then $\mu^*(v) = 2d(v)-6 -n_2(v) - d(v) \cdot \frac{1}{2}  \ge   \frac{3d(v)}{2}-6 -n_2(v) \ge 0$ by (R\ref{g5c7rule:2vx:new'}) and  (R\ref{g5c7rule:good:new'}).
If $v$ is a $5$-vertex with two $2$-neighbors, then $|F^*_5(v)|\le 4$, so  $\mu^*(v) = 4 -2\cdot1 - |F^*_5(v)| \cdot \frac{1}{2} \ge 0$ by (R\ref{g5c7rule:2vx:new'}) and  (R\ref{g5c7rule:good:new'}).
Suppose that $v$ is a $6^+$-vertex. 
We first send charge $1$ from $v$ to each incident edge $vw$, then $vw$ sends either charge $1$ to $w$ if $w$ is a $2$-vertex or charge $\frac{1}{2}$ to each face incident with $vw$ if $w$ is a $3$-vertex.
In this way, we see that $v$ sends charge at most $d(v)$ by applying (R\ref{g5c7rule:2vx:new'}) and (R\ref{g5c7rule:good:new'}).
Thus, if $6\le d(v)\le 8$, then the final charge of $v$ is $\mu^*(v)\ge 2d(v)-6-d(v)\ge 0$.
If $v$ is a $9^+$-vertex, then $v$ sends charge $\frac{1}{3}$ to each incident face by (R\ref{g5c7rule:9vx:new'}) in addition to the charge sent to its incident edges, so $\mu^*(v)=(2d(v)-6)-d(v)-d(v)\cdot \frac{1}{3}\ge 0$.
\end{proof}

\begin{claim}
Each face $f$ has non-negative final charge.
\end{claim}
\begin{proof}
Each face only receives charge, so $\mu^*(f)\geq\mu(f)=d(f)-6\ge 0$ for a $6^+$-face $f$. 

Let $f$ be a $5$-face, which has at least three $4^+$-vertices since no $3^-$-vertices are adjacent to each other. 
If $f$ is incident with at least four $4^+$-vertices, then each of those $4^+$-vertices are giving on $f$, so each of them sends charge at least $\frac{1}{4}$ to $f$ by (R\ref{g5c7rule:good:new'}) and  (R\ref{g5c7rule:bad:new'}), hence $\mu^*(f)\ge -1 + 4\cdot \frac{1}{4}= 0$.
Thus, $f$ has exactly three $4^+$-vertices $u,v,w$. 
Moreover, $f$ has at least two giving $4^+$-vertices $v$ and $w$, so  $f\in F^*_5(v)\cap F^*_5(w)$.
If $v$ is a $9^+$-vertex, then $v$ sends charge at least $\frac{1}{2}+\frac{1}{3}$ to $f$ by (R\ref{g5c7rule:good:new'}) and (R\ref{g5c7rule:9vx:new'}) and $w$ sends charge at least $\frac{1}{4}$ to $f$ by (R\ref{g5c7rule:good:new'}) and (R\ref{g5c7rule:bad:new'}), so $\mu^*(f)\ge -1 + \frac{5}{6}+\frac{1}{4}>0$.
If both $v$ and $w$ are good vertices, then each of $v$ and $w$ sends charge at least $\frac{1}{2}$ to $f$ by (R\ref{g5c7rule:good:new'}), so $\mu^*(f)\ge -1 + 2\cdot \frac{1}{2}= 0$.
So, in the following, we may assume that all giving  vertices on $f$ are $8^-$-vertices, and $f$ has at most one giving good vertex.
So far, we may assume $f$ is a $5$-face $uxvwyu$, where $u$ is a $4^+$-vertex, $v$ is a bad vertex, $d(v),d(w)\le 8$, and $d(x),d(y)\le 3$.

If $f$ is terrible, then by (R\ref{g5c7rule:good:new'}) and  (R\ref{g5c7rule:bad:new'}), each of $v$ and $w$ sends charge at least $\frac{1}{2}$ to $f$. Thus $\mu^*(f)\ge -1 +2\cdot \frac{1}{2}= 0$.

Suppose that $f$ is not terrible. 
If $w$ is a good vertex, then it sends charge at least $\frac{1}{2}$ to $f$ by (R\ref{g5c7rule:good:new'}), so $v$ and $w$ together send charge at least $\frac{1}{2}+\frac{1}{4}$ to $f$. 
If both $v$ and $w$ are bad vertices, then both $v$ and $w$ are $5$-vertices with exactly two $2$-neighbors, since a bad vertex is adjacent to neither a $4$-vertex with a $2$-neighbor nor a $5$-vertex with three $2$-neighbors by Lemma~\ref{lem:combine}.
By (R\ref{g5c7rule:bada:new'}), $v$ and $w$ together send charge at least $\frac{1}{3}+\frac{1}{3}$ to $f$.
Hence, in any case, $v$ and $w$ together send charge at least $\frac{2}{3}$ to $f$.

We will show that $u$ always sends charge at least $\frac{1}{3}$ to $f$, so  the final charge of $f$ is non-negative. 
If $u$ is a $9^+$-vertex, then $u$ sends charge $\frac{1}{3}$ to $f$ by (R\ref{g5c7rule:9vx:new'}).
If $u$ is an $8^-$-vertex, then $x$ or $y$ is a $3$-vertex, since $f$ is not terrible and $d(v), d(w)\leq 8$, so $u$ is a giving vertex on $f$. 
Thus, if $u$ is a good vertex, then $u$ sends charge $\frac{1}{3}$ to $f$ by (R\ref{g5c7rule:good:new'}).
If $u$ is a bad vertex, then $u$ has both a $2$-neighbor and a $3$-neighbor, so 
Lemma~\ref{lem:2neighbor-new} implies $n_2(u)+n_3(u) \le 2d(u)-7$.
Hence, $u$ is a $5$-vertex with exactly two $2$-neighbors, so $u$ sends charge at least $\frac{1}{3}$ to $f$ by (R\ref{g5c7rule:bada:new'}).
\end{proof}

\section{Concluding remarks}

We remark that most of our techniques and ideas come from our own paper~\cite{arXiv_ChChKwBo} on odd coloring, but we had to enhance them significantly for proper conflict-free coloring. 
The analysis is more complicated and requires extra effort since proper conflict-free coloring is a strengthening of odd coloring. 
An {\it odd coloring} of a graph is a proper coloring such that each non-isolated vertex has a color appearing an odd number of times on its neighborhood.
We find it quite interesting that the threshold on the maximum average degree guaranteeing an odd $c$-coloring is exactly the same as the maximum average degree threshold guaranteeing a proper conflict-free $c$-coloring. 
However, the two parameters indeed have a difference when applied to planar graphs. 
As mentioned in the introduction, Fabrici et al.~\cite{arXiv_FaLuRiSo} constructed a planar graph that does not have a proper conflict-free 5-coloring, whereas it is conjectured that planar graphs have an odd 5-coloring~\cite{petrusevski2021colorings}. 

In~\cite{arXiv_FaLuRiSo}, it is conjectured that every planar graph has a proper conflict-free $6$-coloring, and it is proved that a planar graph has a proper conflict-free 8-coloring.
Let $g_c$ be the minimum $k$ such that all planar graphs with girth at least $k$ has a proper conflict-free $c$-coloring. 
Theorem~\ref{thm:mad:c6} implies $g_6\leq 8$ and Theorem~\ref{thm:g5c7} is equivalent to  $g_7\leq 5$.
It would certainly be interesting to determine the exact value of $g_c$ for each $c\in\{4,5, 6,7\}$.
Note that $g_c=3$ for all $c\geq 8$ since it is known that every planar graph has a proper conflict-free 8-coloring and no such value exists for $c\leq 3$ since every cycle on $3k\pm 1$ vertices does not have a proper conflict-free 3-coloring. 

As pointed out in~\cite{arXiv_CaPeSk}, if one can prove $g_4\leq 6$, then it would be a strengthening of the Four Color Theorem~\cite{1977ApHa,1977ApHaKo,1997RoSaSeTh}. 
Note that our Theorem~\ref{thm:mad:c4} implies $g_4\leq 10$.

\section*{Acknowledgements}
Eun-Kyung Cho was supported by Basic Science Research Program through the National Research Foundation of Korea (NRF) funded by the Ministry of Education (NRF-2020R1I1A1A0105858711).
Ilkyoo Choi was supported by the Basic Science Research Program through the National Research Foundation of Korea (NRF) funded by the Ministry of Education (NRF-2018R1D1A1B07043049), and also by the Hankuk University of Foreign Studies Research Fund.
Boram Park was supported by Basic Science Research Program through the National Research Foundation of Korea (NRF) funded by the Ministry of Science, ICT and Future Planning (NRF-2018R1C1B6003577).

\bibliographystyle{abbrv}
\bibliography{ref}

\end{document}